\documentclass{amsart}

\usepackage{amssymb, amsmath}

\begin{document}

\def\frg{\frak g}
\def\frU{\frak U}
\def\frM{\frak M}
\def\deg{\mathrm {deg}}
\def\wh{\widehat}
\def\Ind{\mathrm {Ind}}
\def\wt{\widetilde}
\def\ov{\overline}
\def\wt{\widehat}

\def\SU{\mathrm {SU}}
\def\SL{\mathrm {SL}}
\def\PU{\mathrm {PU}}
\def\PSL{\mathrm {PSL}}
\def\U{\mathrm  U}
\def\OO{\mathrm  O}
\def\Sp{\mathrm  {Sp}}
\def\SO{\mathrm  {SO}}
\def\SOS{\mathrm {SO}^*}

\def\bs{\bigstar}

\def\C{\mathbb C}
\def\R{\mathbb R}

\def\epsilon{\varepsilon}

\renewcommand{\Re}{\mathop{\rm Re}\nolimits}

\renewcommand{\Im}{\mathop{\rm Im}\nolimits}

\def\ch{\cosh}
\def\sh{\sinh}

\begin{center}
{\bf\large
$K$-finite matrix elements of 
irreducible Harish--Chandra modules
are hypergeometric}

\sc\large
Neretin Yu.A.\footnote{supported by grant NWO--047.017.015 }  
\end{center}

{\small We show that each $K$-finite matrix
element of an irreducible Harish-Chandra module
can be obtained from spherical functions
by a finite collection of operations.}

\medskip

{\bf 1. Notation.}
Let $G$ be a linear semisimple Lie group,
let $K$ be the maximal compact subgroup.
Let $\frg$ be the Lie algebra of $G$,
let $\frU(\frg)$ be its universal
enveloping algebra.
%For $p\in \frU(\frg)$ denote by $\deg(\cdot)$
%its  degree.
Let $L_X$, $R_X$, where $X\in\frg$,
be the left and right Lie derivatives on $G$.
Denote by $\frU_l(\frg)$ (resp., $\frU_r(\frg)$)
the  algebra of differential operators on $G$
generated by  left  (resp., right) derivatives.

By $\Psi_s(g)$ we denote the spherical functions on $G$,
$s$ is the standard parameter of a spherical
function
(see \cite{Hel}).
 
For a finite-dimensional representation
$\xi$ of $G$, denote by $\frM(\xi)$ the space 
of finite linear combinations of 
 matrix elements of $\xi$.

\smallskip

{\bf 2. Formulation of the result.}
Let $V$ be an irreducible Harish-Chandra module
(see, for instance, \cite{Kna})
over $G$, denote by $\pi(g)$ operators of representation
of $G$ in some completion of $V$.
Let $\sigma$ ranges in the set $\wh K$
of all irreducible representations
of $K$.
Let $V=\oplus_\sigma V_\sigma$ be the decomposition
of $V$ into a direct sum of $K$-isotipical components.

\smallskip

{\sc Proposition}.  a)  {\it Let $V$ be an irreducible
Harish-Chandra module in a general position.
Let $V^\circ$ be the dual module.
Let $v\in V_\sigma$, $w\in V^\circ_\tau$.
There exists an irreducible finite dimensional representation
$\xi$ of $G$ and $s$ such that the matrix element
$\{\pi(g)v,w\}$ is a finite sum
$$
\{\pi(gv),w\}=\sum_j h_j(g)\cdot p_j q_j \Psi_s(g)
$$
where $h_j\in \frM(\xi)$, $p_j\in\frU_l(\frg)$, $q_j\in\frU_r(\frg)$.
}

b) {\it  For an arbitrary Harish-Chandra module,
each $K$-finite matrix element admits a representation
$$
\{\pi(g)v,w\}=\lim_{s\to s_0}
\sum_j h_j(s;g)\cdot p_j(s) q_j(s) \Psi_s(g)
,$$
where $h_j(s;g)$ is an element $\frM(\xi)$
  depending in a parameter $s$, and $p_j(s)\in\frU_l(\frg)$, $q_j(s)\in\frU_r(\frg)$.
Moreover, the degrees of $p_j(s)$, $q_j(s)$,
and the number of summands are uniformly
bounded in $s$ (for  fixed $V$, $\sigma$, $\tau$).
} 

%\smallskip

%{\sc Remark.} The construction of these expressions
%are quasi-explicit (see the proof), but of course they are not
%really explicit formulae. Nevertheless, these expressions
%allow to obtain a priory statements on matrix coefficients
%(for instance, on domains of holomorphy,
%see, \cite{AG}, \cite{KG}).

\smallskip

{\bf 3. Proof.}
First, we introduce additional notation.

 Denote by $P$ the minimal parabolic subgroup
in $G$, consider the decomposition 
$P=MAN$, where $N$ is the  nilpotent radical,
$MA$ is the Levi factor of $P$,
$M$ is the  compact subgroup, $A\simeq(\R_+^*)^k$ 
is  the vector subgroup.

 Consider the flag space
$G/P$, consider the corresponding
Grassmannians, i.e., factor-spaces $G/Q_\alpha$,
where $Q_\alpha\supset P$ are maximal parabolics in $G$.
Equip all the spaces $G/P$, $G/Q_\alpha$
with $K$-invariant measures (this  allows
to define Jacobians below).
For $\omega\in G/P$, denote by $\omega_\alpha$
its image under the map $G/P\to G/Q_\alpha$.

For  $g\in G$, denote by $J_\alpha(g,\omega)$ the Jacobian
of the transformation $g:G/Q_\alpha\to G/ Q_\alpha$ at the
point $\omega_\alpha$.
  
Let $\mu$ be a character $A\to \C^*$, let 
$\tau$ be an irreducible representation of $M$.
We denote by $\mu\otimes \tau$ the representation
of $P=MAN$, that is $\mu$ on $A$, $\tau$ on $M$
and is trivial on $N$. By $\Ind_P^G(\mu\otimes\tau)$
we denote the representation of $G$ induced from
$\mu\otimes\tau$, i.e., a representation of principal
(nonunitary) nondegenerate series.

By  
 the Subquotient Theorem,
each Harish-Chandra module can be realized
as a subquotient (and even as a subrepresentation)
in some representation $\Ind_P^G(\mu\otimes\tau)$,
see, for instance,  \cite{Kna}.
Hence, it is sufficient to prove the statement
on matrix elements for the representations
$\Ind_P^G(\mu\otimes\tau)$ (they can be reducible).

$A^\bs$.
Let $\rho$ be a spherical representation with parameter
$s$. Let $h\in V$ be a spherical vector, let $h^\circ$
be the spherical vector in   $V^\circ$.
Vectors $v\in V_\sigma$, $w\in (V^\circ)_\tau$ can be represented
in the form
$$
v=\Bigl( \sum a_\alpha \prod_j X_{\alpha j}\Bigr)\cdot h ,\qquad
 w=\Bigl(\sum b_\beta \prod_i Y_{\beta i}\Bigr)\cdot h^\circ
$$
for some $X_{\alpha j}$, $Y_{\beta i}\in\frg$.
Thus
\begin{multline*}
\Bigl\{ \pi(g)v,w\Bigr\}=
\Bigl\{\pi(g)
\Bigl( \sum a_\alpha \prod_j X_{\alpha j}\Bigr)\cdot h,
\Bigl(\sum b_\beta \prod_i Y_{\beta i}\Bigr)\cdot h^\circ
\Bigr\}
=
\\
=
\Bigl\{
\Bigl(\sum b_\beta \prod_i (- Y_{\beta i}) \Bigr)
\pi(g)\Bigl( \sum a_\alpha \prod_j X_{\alpha j}\Bigr)\cdot h,
 h^\circ
\Bigr\}
=\\=
\Bigl(\sum b_\beta \prod_i (-L_{ Y_{\beta i} }) \Bigr)
\Bigl( \sum a_\alpha \prod_j (R_{X_{\alpha j}}\Bigr)
\Psi_s(g)
\end{multline*}

$B^\bs$.
Consider an
induced representation
$
\pi
=\Ind_P^G(\chi\otimes 1)
,$
where $1$ denotes one-dimensional representation
of $M$.
If $\chi=\chi_s$ is 
in a general position (in fact $s\in\C^k$ is outside
a locally finite family of complex hyperplanes), 
then $\pi$ is an irreducible
 spherical representation.
This situation was considered
in $A^\bs$.
Now examine the case of reducible  $\pi$.
For this, we must follow continuity
of matrix elements as functions of
parameters $s$.

For $t_\alpha\in\C$ define the representation
$$\rho_t(g)f(\omega):= f(g\omega)
 \prod J_\alpha(g,\omega)^{t_\alpha} $$
 of $G$ in the space of
functions on $G/P$.
It can be readily checked that this family
of representations coincides with the family
$\Ind_P^G(\chi\otimes 1)$, where $\chi=\chi_s$ 
ranges in all the characters of $A$
(and the dependence $s=s(t)$ 
is some linear transformation).

Thus, we obtain a realization of the family
$\Ind_P^G(\chi\otimes 1)$ such that the action of $K$
is independent in $\chi$ and operators of representation
are continuous functions in $s$.

$C^\bs$.
Let $\xi$ be an irreducible
finite-dimensional representation
of $G$ in the space $H$.
Following  \cite{Kos},
we consider
the tensor product
\begin{equation}
\pi\otimes\xi=
 \Ind_P^G(\chi) \otimes\xi
=\Ind_P^G(\chi\otimes \xi\Bigr|_P)
\end{equation}
The representation 
$ \xi\Bigr|_P$
is reducible and it admits
a finite  filtration with irreducible
subquotients
$$H_1\supset H_2\supset H_3\supset\dots$$
The nilpotent subgroup $N\subset P$ acts in the subquotients
$H_j/H_{j+1}$ in a trivial way. The
representations of subgroup $MA\subset P$ 
in $H_j/H_{j+1}$
have the form
$\mu_j\otimes\tau_j$
for some characters $\mu_j$ of $A$ and some 
irreducible representations $\tau_j$ of $M$. 

 Thus, the representation $\pi\otimes \rho$ has a filtration,
whose subquotients
are representations of principal series having the form
$
\Ind_P^G([\chi\cdot\mu]\otimes \tau)
$.

$D^\bs$. Fix a representation 
$\wt\tau$ of $M$ and a character $\wt\chi$ of $A$.
We intend to realize $\Ind_P^G(\wt\mu\otimes\wt \tau)$ as subquotient
in an appropriate tensor
product (1).

We can choose a  representation
$\xi$ of $G$ 
such that the restriction of $\xi$ to $M$ contains $\wt\tau$.%
\footnote{A proof. Denote by $G_c$ the compact form of $G$. 
Consider the induced representation $\Ind_L^{G_c}(\wt\tau)$.
Let $\xi$ be its irreducible subrepresentation.
We consider $\xi$ as a representation of $G$.}
Then restriction of $\xi$ to $P=MAN$ contains a subquotient
of form $ \wt\mu\otimes\wt\tau$
with a certain character $\mu$.

Next, we choose
a character $\chi$ of $A$
 such that $\chi\cdot \wt\mu=\wt\chi$.
Thus we obtain that 
$\Ind_P^G(\chi\otimes 1)\otimes \xi$ contains a given representation
$\Ind_P^G(\wt\mu\otimes \wt\tau)$ as a subquotient.

$E^\bs$. $K$-finite matrix elements of $\Ind_P^G(\wt\mu\otimes \wt\tau)$ 
are contained in $K$-finite matrix elements of 
$\Ind_P^G(\chi)\otimes \xi$. The latter  matrix elements
are finite linear combinations  
of products of $K$-finite matrix elements of $\Ind_P^G(\chi)$
and matrix elements of $\xi$. 
This finishes proof of a) and b).

\smallskip

{\bf 4. An application. Domains of holomorphy of matrix coefficients.}
%Our Proposition  allows to produce a priory statements
%on behavior of matrix elements, since spherical
%functions are relatively
% well-investigated (see, for instance, \cite{HO},
% \cite{HS},
%Chapter 1).

\smallskip

{\sc Corollary.} {\it Each domain of holomorphy
$\Omega\subset G_\C$
of all the spherical functions
is a domain of holomorphy of all the
$K$-finite
matrix elements of all the irreducible
Harish-Chandra modules over $G$.
}

\smallskip

In particular, all such matrix elements are
holomorphic in the Akhiezer--Gindikin domain 
 \cite{AG} (this is obtained in \cite{KS}).

{\sc Corollary.} {\it There is a submanifold
$Y\subset G_\C$, such that for each irreducible
Harish-Chandra module over $G$, the matrix-valued
function $g\mapsto \rho(g)$ can be extended to a holomorphic function
on the universal covering space of $G_\C\setminus Y$.}

{\sc Proof.}  
Spherical functions on $G$ are multivalued
holomorphic functions
on $G_\C$ having singularities (branching)
on a prescribed manifold $Y\in G_\C$,
see \cite{HO}.
Hence, for
Harish-Chandra modules  in a general position,
there is nothing to prove.

To prove the  statement for exceptional values of 
$s$, we must follow details of
$B^\bs$.

Denote by $V$ the space of $K$-finite functions
on the flag space $G/P$ (see $B^\bs$).
Denote by $1\in V$ the function $f(\omega)=1$.
Fix $\sigma, \tau\in \wh K$.
Fix $w\in V_\sigma$, $w^\circ\in V_\tau$.
Let consider the representations
$\Ind_P^G(\chi_s\otimes 1)$ and
 follow a behavior
of the corresponding matrix element
as a function of $s$.

 Denote by $\frU^N(\frg)\subset \frU(\frg)$
the subspace consisting of all elements of
degree $\leqslant N$. 
Let $N$ be sufficiently large,
such that $\frU^N(\frg)\cdot 1$
contains the whole  subspace 
$V_\sigma$
for all generic characters $\chi$.
Consider a collection $r_1$, 
$r_2, \dots\in\frU^N$ such that
for some generic $\chi$

1. $r_j \cdot 1$ are linear independent in 
$\Ind_P^G(\chi_s)$

2. their linear combinations contains
$V_\sigma$.

These properties remain valid
for all $s$ outside
a certain algebraic submanifold
$\mathcal M$ in the space of parameters.
Thus we express our vector $w$
 as a linear combination of $r_j$,
$w=\sum c_j(s)r_j(s)$,
where $c_j$ are certain rational functions.

Applying $A^\bs$, we obtain, that
our matrix element
 has a form
$%\begin{equation}
\Xi(s) \Psi_s(g)
$, %,\end{equation}
where $\Xi(s)$ is an element of 
$\frU_l\otimes\frU_r$ depending 
rationally in $s$.  

Now let $s_0$ be an exceptional value of
$s$. Let $w\in V_\sigma$, $w^\circ\in V_\tau$.
Consider a holomorphic curve
$\gamma(\epsilon)$ (with $\gamma(0)=s_0$) avoiding singularities
of $\Xi$ and singular values of the parameter
$s$.
A priory, the function 
$$
F(\epsilon,g)
=\Xi(\gamma(\epsilon))
\Psi_{\gamma(\epsilon)}(g)
$$
is holomorphic in the domain
$0<|\epsilon|<\delta$, $g\in G_\C\setminus Y$
and has a pole on the submanifold $\epsilon=0$.
But we know, that that $F(\epsilon,g)$ has 
a finite limit as $g\in G$
 and $\epsilon\to 0$.
Hence it has no pole, and hence it is holomorphic
at $\epsilon=0$.
In particular it is holomorphic on the submanifold
$\epsilon=0$, and this is the desired statement.
\hfill$\square$

A product of such matrices $\rho(g_1)\cdot\rho(g_2)$
generally is divergent but sometimes
it is well-defined(see \cite{Nel}, \cite{KS}).
For each $X\in\frg_\C$, $g\in G_\C\setminus Y$, we have
$$
\frac{d}{d\epsilon}\rho(\exp(\epsilon X)g)
=\rho(X)\rho(g)
$$

{\sc Remark.}
 There are  exceptional situations, when
a unitary  representation
admits a continuation to the whole complex
group or its subsemigroup,
apparently these cases are well-understood,
see \cite{Ols}, \cite{Ner1}, 
\cite{Ner2}, 
(Sections 1.1, 4.4, 5.4, 7.4-7.6, 9.7),
\cite{Lit}, \cite{Goo}.
May be there are other (non-semigroup) cases of 
unexpectedly large
(non-semigroup) domain of
holomorphy. As far as I know, this 
problem never was considered.

\smallskip 

{\bf 5. Nonlinear semisimple Lie groups.}
For universal coverings of the groups
$\SU(p,q)$, $\Sp(2n,\R)$, $\SOS(2n)$ our construction
survives, we only must replace spherical functions
by appropriate Heckman--Opdam hypergeometric
functions, see \cite{HS}, Chapter 1.

I do not know, is it possible to
express matrix elements of 
universal covering groups of $\SO(p,q)$ and $\SL(n,\R)$
in the terms of Heckman--Opdam hypergeometric functions.

{\tt Math.Physical Group, ITEP, B.Cheremushkinskaya, 25,
Moscow, 117259, Russia

Fakult\"at f\"ur Mathematik, Universit\"at Wien,
Nordbergstrasse~15, A-1090 Wien, Austria}

e-mail: {\tt neretin@mccme.ru}

\end{document}